\documentclass[12pt]{article}
\usepackage{amsmath,amssymb}
\usepackage{authblk}
\usepackage{color}
\usepackage{mathrsfs}
\usepackage{mathptmx}
\usepackage{verbatim}
\usepackage{epsfig}
\usepackage{CJK}
\usepackage{bm}
\usepackage{amsfonts}
\usepackage{amsthm}
\input{epsf.sty}
\usepackage{epsfig}
\usepackage{setspace}

\def\<{\langle}
\def\>{\rangle}

\def\be{\begin{equation}}
\def\ee{\end{equation}}
\def\ba{\begin{array}}
\def\ea{\end{array}}

\newtheorem{theorem}{Theorem}[section]

\newtheorem{lemma}{Lemma}[section]
\newtheorem{proposition}{Proposition}[section]

\theoremstyle{definition}
\newtheorem{remark}{Remark}[section]

\newtheorem{definition}{Definition}[section]

\numberwithin{equation}{section}

\usepackage{amsmath}
\usepackage{amsthm}
\usepackage{amsfonts}
\usepackage{pifont}
\usepackage{color}
\usepackage{bbm}
\usepackage{CJK}
\usepackage{mathrsfs}
\usepackage{fancyhdr}
\usepackage{graphicx}
\usepackage{pdfsync}
\usepackage{psfrag}
\usepackage{time}
\usepackage{txfonts}

\usepackage{stmaryrd}
\usepackage{mathrsfs}
\usepackage{txfonts}
\usepackage{amsfonts}

\usepackage{indentfirst}

\def\be{\begin{equation}}
\def\ee{\end{equation}}
\def\br{\begin{eqnarray}}
\def\er{\end{eqnarray}}

\title{Weak KAM Theorem for a Class of Infinite-Dimensional Lagrangian systems
\thanks
{This work was supported by the NNSF (11371132) and Key Laboratory of High Performance Computing and Stochastic Information Processing }}

\author[1]{Guanghua Shi\thanks{12110180067@fudan.edu.cn}}
\author[2]{Cheng Yang\thanks{11110180006@fudan.edu.cn}}
\affil[1]{College  of Mathematics and Computer  Science, Hunan Normal University, \newline Hunan  410006, People's Republic of China}
\affil[2]{School of Mathematical Sciences, Fudan University, Shanghai 200433,  \newline People's Republic of China}

\begin{document}

\maketitle

\date{}

\begin{quote}
\small {\bf Abstract.}  In this paper, we study infinite-dimensional Lagrangian systems where the potential functions are periodic, rearrangement invariant and weakly upper semicontinuous. And we prove  that there exists a calibrated curve for every $M\in L^{2}(I)$.

\end{quote}

\tableofcontents

\section{Introduction}
In the nineties of last century, Fathi established the weak KAM theory for autonomous Tonelli Lagrangian systems on finite-dimensional compact manifolds (see \cite{fat1},\cite{fat2},\cite{fat3}). He introduced the weak KAM solution and found its relation with the viscosity of Hamilton-Jacobi equation \cite{cra}. To prove the existence of the weak KAM solution, Fathi adopted the method of Lax-Oleinik semigroup which is well-known in PDE and in Calculus of Variations. Specifically, one need to prove that the semigroup is convergent with respect to certain topology. Since then many people have studied the weak KAM theory for more general Lagrangian systems with relaxed conditions.

Firstly, some scholars considered the Lagrangian systems on noncompact manifold, refer to \cite{con}, \cite{fat4}. In \cite{fat4}, Fathi and Maderna proved that the Lax-Oleinik semigroup is convergent with respect to the compact open topology. While, in \cite{con}, a weak KAM solution was directly constructed by using semistatic curve. When the Lagrangian is nonautonomous, the Lax-Oleinik semigroup fails to converge, refer to \cite{fat6}. To overcome this difficulty, Wang and Yan \cite{wan} defined a new kind of semigroup which is convergent. In addition, Maderna \cite{mad} studied the weak KAM theory for $N-$body problems, in which the Lagrangians have singularities, not satisfy the Tonelli conditions. To deal with the singularities, he used cluster partitions techniques and obtained precise estimates of the action functional.

The generalized weak KAM results mentioned above are on finite-dimensional manifolds. As to the infinite-dimensional cases, Gangbo \cite{gan1}, \cite{gan2} made a remarkable work.
In the Hilbert space $L^{2}(I)$ ($I:=(0,1)$), Gangbo \cite{gan1} considered the Lagrangian function
\begin{equation*}
  L(M,N)=\frac{1}{2}||N||^{2}_{L^{2}(I)}-\mathcal{W}(M), \quad (M,N)\in L^{2}(I)\times L^{2}(I),
\end{equation*}
where the potential $\mathcal{W}$ is continuous, periodic and rearrangement invariant (these concepts will be explained in Section 2). Due to Remark \ref{r22}, the Lagrangian sysem can be seen as a system on the Wasserstein space $\mathcal{P}(\mathbb{T}^{1})$, which is compact. Then he made use of the approximation method to obtain a solution $u$ satisfying
\begin{itemize}
  \item $u$ is dominated by $L+\lambda$ on $L^{2}(I)$;
  \item If $M\in L^{2}(I)$ is monotone nondecreasing, there exists a calibrated curve $\sigma^{M}:(-\infty,0]$ $\rightarrow L^{2}(I)$ with $\sigma^{M}(0)=M$.
\end{itemize}
Besides, he applied this result to studying the dynamic behavior of Vlasov system. By the definition of the weak KAM solution, $u$ is not a weak KAM solution on the whole Hilbert space. So there is a natural question whether one can find a calibrated curve for every $M\in L^{2}(I)$. To conclude it, the key step is to check that the Tonelli Theorem holds true in $L^{2}(I)$, not a subset of $L^{2}(I)$ in \cite{gan1}. 
\par
In this paper, we will give a positive answer (see Theorem \ref{t64}) to this question by assuming an extra assumption on the potential functions:
\begin{itemize}
  \item $\mathcal{W}$ is weakly upper semicontinuous.
\end{itemize}
Note that the weakly upper semicontinuity of potential plays an important role in the existence of minimizers for the action functional $\int L ds$.

Our proof is based on the method of Lax-Oleinik semigroup. The method consists of two steps. The first step is to show that the infimum in (\ref{e41}) can be realized which is the major difficulty in this paper. For this, we introduce an analogous priori compactness (refer to Lemma \ref{l32}). The second step is to find the fixed point of Lax-Oleinik semigroup. This step is similar to that in \cite{fat4}.

Finally, we give an outline of the remaining parts. In Section 2, some concepts and known facts are displayed. Then we generalize the Tonelli Theorem in finite space to the Hilbert space in Section 3. The following section is devoted to the Lax-Oleinik semigroup. In Section 5, we show some properties of the dominated function. In the last section, we prove the weak KAM Theorem.

\section{ Preliminaries and some known results  }
At first, we introduce the concepts of periodic and rearrangement invariant. In the finite dimensional space $\mathbb{R}^{n}$, a function $f$ is said to be periodic if it satisfies the following condition
\begin{equation*}
  f(x+k)=f(x)\quad \forall x\in \mathbb{R}^{n},k\in\mathbb{Z}^{n}.
\end{equation*}
To define the periodic function in the Hilbert space, we consider the group below
$$L_{\mathbb{Z}}^{2}(I):=\{M\in L^{2}(I)\;|\;M(x)\in\mathbb{Z},\;\forall x\in I\}$$
which can be seen as the extension of the group $\mathbb{Z}^{n}$.
\begin{definition}(\textbf{Periodic})
We say that a function $u:L^{2}(I)\rightarrow\mathbb{R}$ is periodic, if
\begin{equation*}
  u(M+N)=u(M)\quad \forall M\in L^{2}(I),N\in L^{2}_{\mathbb{Z}}(I).
\end{equation*}
\end{definition}

Let $P_{n}$ denote the set of permutation of $n$ letter. For $x=(x_{1},\cdots,x_{n})\in\mathbb{R}^{n}$ and $\sigma\in P_{n}$, $\sigma (x)$ stands for the vector obtained by permuting the components of $x$ according to $\sigma$. If the function $f$ defined in $\mathbb{R}^{n}$ satisfies
$$f(x)=f(\sigma(x))\quad \forall x\in \mathbb{R}^{n}, \sigma\in P_{n},$$
we say that $f$ is rearrangement invariant. To extend the concept of rearrangement invariant function to Hilbert space, we firstly generalized the group of permutation of $n$ letter. Let $\mathcal{G}$ denote the collection of bijections $G:[0,1]\rightarrow[0,1]$ such that $G,G^{-1}$ are Borel and push $\nu_{0}$ forward to itself, where $\nu_{0}$ is the one-dimensional Lebesgue measure restricted to $I$.
\begin{definition}(\textbf{Rearrangement invariant})
We say that a function $u:L^{2}(I)\rightarrow\mathbb{R}$ is rearrangement invariant, if
\begin{equation*}
  u(M\circ G)=u(M)\quad \forall M\in L^{2}(I), G\in \mathcal{G}.
\end{equation*}
\end{definition}

In the following part, we introduce the Wasserstein space $\mathcal{P}(\mathbb{T}^{1})$. The element in this space is Borel probability measure on the torus $\mathbb{T}^{1}$. Given $\mu,\nu\in\mathcal{P}(\mathbb{T}^{1})$, the distance between $\mu$ and $\nu$ is
$$\inf_{\gamma\in\Gamma(\mu,\nu)}\int_{\mathbb{T}^{1}\times\mathbb{T}^{1}}|x-y|^{2}_{\mathbb{T}^{1}}d\gamma(x,y),$$
where $|\cdot|_{\mathbb{T}^{1}}$ is the distance on the torus $\mathbb{T}^{1}$ and $\Gamma(\mu,\nu)$ is the set of Borel measures on $\mathbb{T}^{1}\times\mathbb{T}^{1}$ which have $\mu$ and $\nu$ as marginal.
\begin{remark}\label{r21}
 $\mathcal{P}(\mathbb{T}^{1})$ is a compact, complete, separable metric space, refer to \cite{gan1}.
\end{remark}
\begin{remark}\label{r22}
If a continuous function in $L^{2}(I)$ is periodic and rearrangement invariant, it is well-defined on $\mathcal{P}(\mathbb{T}^{1})$. For details, you can refer to Proposition 2.8 and Lemma 2.14 in \cite{gan1}.
\end{remark}
According to the remarks above, there exists a constant $K_{0}$ such that $$|\mathcal{W}(M)|\leq K_{0},\;\forall M\in L^{2}(I).$$

At last, we introduce 2-absolutely continuous curve \cite{amb}. Suppose $(\mathcal{S},dist)$ is a complete metric space and a curve $\sigma:[0,T]\ni t\mapsto\sigma_{t}\in\mathcal{S}$. If there exists $\beta\in L^{2}[0,T]$ such that dist$(\sigma_{t},\sigma_{s})\leq\int_{s}^{t}\beta(u)du$ for every $s<t$ in $[0,T],$ we say that $\sigma$ is 2-absolutely continuous. We denote by $AC^{2}(0,T;\mathcal{S})$ the set of 2-absolutely continuous curves. By H$\ddot{\text{o}}$lder inequality, we know that if $\sigma$ is 2-absolutely continuous, it is $\frac{1}{2}-$H$\ddot{\text{o}}$lder continuous.

\section{Tonelli Theorem }
In this section, we show a generalized version of Tonelli Theorem in Hilbert space. To prove it, we give a useful lemma which is about compactness of a sequence of curves in $AC^{2}(0,T;L^{2}(I))$. Firstly, one defines the topology $\tau$ on $AC^{2}(0,T;L^{2}(I))$.
\begin{definition}   \label{d31}
 Let $\sigma^{k},\sigma\in AC^{2}(0,T;L^{2}(I)),$ we say that $\sigma^{k}$ converges to $\sigma$ with respect to the topology $\tau$, denoted by $\sigma^{k}\stackrel{\tau}{\rightarrow}\sigma,\;\;if\;\;\sigma^{k}(t)\rightharpoonup\sigma(t)\;\;in\; L^{2}(I),\;\forall t \in[0,T].$
\end{definition}

As we know that to prove that the set of curves in finite manifold is compact, we often use Ascoli-Arzal\`{a} Theorem. Fortunately, there is a refined version of Ascoli-Arzal\`{a} Theorem in a metric space.
\begin{lemma} \label{l31}
 Let $T>0$, and $K\subset\varphi$ be weak compact, where $\varphi$ is a metric space and the distance is lower semicotinuous with respect to the weak topology that is Hausdorff. Let $u_{n}:[0,T]\rightarrow\varphi$ be curves satisfying
 \begin{itemize}
\item[a)] $u_{n}(t)\in K \;\;\forall n\geq1,\;\;\forall t\in[0,T],$
\item[b)] $\limsup_{n}d(u_{n}(t),u_{n}(s))\leq\omega(t,s), \;\;where \;\;\lim_{(s,t)\rightarrow(r,r)}\omega(t,s)=0.$
\end{itemize}
 Then there exist a subsequence $u_{n_{k}}$ and a limit curve $u:[0,T]\rightarrow\varphi$\;such that $u_{n_{k}}(t)\rightarrow u(t)$ $\forall t\in[0,T]$ with respect to the weak topology\;and\;$u:[0,T]\rightarrow\varphi$\;is strong continuous.
\end{lemma}
The lemma can be found in page 69 of \cite{amb}.
\par
According to the above lemma, we can prove the useful lemma.
\begin{lemma} (\textbf{A Priori Compactness})  \label{l32}
 Let $\sigma^{k}:[0,T]\rightarrow L^{2}(I)$ be a sequence of 2-absolutely continuous curves. We suppose that the sequence $\sigma^{k}(t_{0})$ is bounded in norm for some $t_{0}\in [0,T]$. If the curves satisfy $\sup_{k}\int_{0}^{T}\parallel \dot{\sigma}^{k}(t)\parallel_{L^{2}(I)}^{2}dt<+\infty$, then there exist $\sigma\in AC^{2}(0,T;L^{2}(I))$ and a subsequence still denoted by $\{\sigma^{k}\}$ such that $\sigma^{k}\stackrel{\tau}{\rightarrow}\sigma.$
\end{lemma}
\noindent\textbf{Proof.}
Set \;$K_{1}^{2}:=\sup_{k}\int_{0}^{T}\parallel \dot{\sigma}^{k}(t)\parallel_{L^{2}(I)}^{2}dt.$ For every $s<t$ in $[0,T]$, it is easy to obtain that
\begin{eqnarray*}
  \parallel\sigma^{k}(s)-\sigma^{k}(t)\parallel &\leq& \int_{s}^{t}\parallel \dot{\sigma}^{k}(\tau)\parallel_{L^{2}(I)}d\tau \\
   &\leq& \sqrt{t-s}(\int_{0}^{T}\parallel \dot{\sigma}^{k}(t)\parallel_{L^{2}(I)}^{2}dt)^{\frac{1}{2}} \\
   &\leq&  |M|\sqrt{t-s}\leq |K_{1}|T.
\end{eqnarray*}
According to the above inequality, we not only get $$\limsup_{k}\parallel\sigma^{k}(s)-\sigma^{k}(t)\parallel\leq|M|\sqrt{t-s},$$
but also find a constant $R>0$ such that $\sigma^{k}(t)\in B_{R}$ \footnote{$B_{R}$ denotes the closed ball of radius $R$ and of center at origin in $L^{2}(I)$. }$,\forall t\in[0,T]\text{ and } k\geq1$. We also know that the weak topology of $B_{R}$ is Hausdorff due to the fact that the weak topology of a bounded closed  set in $L^{2}(I)$ can be metrizable. To apply Lemma \ref{l31}, it remains to prove the norm of $L^{2}(I)$ is weak lower semicontinuous with respect to the weak topology in $L^{2}(I)$. This is equivalent to prove that for any $\alpha\geq0$, the set
$$\Sigma_{\alpha}:=\{M\in L^{2}(I):\;\;\parallel M\parallel_{L^{2}(I)}\leq\alpha\}$$ is  weak closed. It is true since $\Sigma_{\alpha}$ is convex and closed.

By Lemma \ref{l31}, there exist a subsequence $\sigma^{k}$\;and a limit curve\;$\sigma$\;such that $\sigma^{k}\stackrel{\tau}{\rightarrow}\sigma$ and $\sigma:[0,1]\rightarrow B_{R}$\;is continuous with respect to $L^{2}$ norm.
Now we consider the time derivative of $\sigma$ in the sense of distribution. For any $\varphi\in C_{c}^{\infty}([0,T]\times I)$, one has
\begin{eqnarray*}
 <\dot{\sigma},\varphi> &=& -<\sigma,\dot{\varphi}>_{L^{2}([0,T]\times I)}\\
 &=& -\int_{0}^{T}<\sigma(t),\dot{\varphi}(t)>_{L^{2}(I)}dt \\
   &=&-\lim_{k}\int_{0}^{T}<\sigma^{k}(t),\dot{\varphi}(t)>_{L^{2}(I)}dt \\
  &=&\lim_{k}\int_{0}^{T}<\dot{\sigma}^{k}(t),\varphi(t)>_{L^{2}(I)}dt \\
   &\leq& \lim_{k}(\int_{0}^{T}\parallel \dot{\sigma}^{k}(t)\parallel_{L^{2}(I)}^{2}dt)^{\frac{1}{2}}\parallel\varphi\parallel_{L^{2}([0,T]\times I)} \\
 &\leq&|K_{1}|\parallel\varphi\parallel_{L^{2}([0,T]\times I)}.
\end{eqnarray*}
So $\dot{\sigma}$ is a bounded linear functional on $L^{2}([0,T]\times I)$. By Riesz Representation Theorem, we know $$\dot{\sigma}\in L^{2}([0,T]\times I)\cong L^{2}([0,T];L^{2}(I)).$$
Finally, one easily gets $$\parallel \dot{\sigma}(t)\parallel_{L^{2}(I)}\in L^{2}(0,T).$$
Since
$$||\sigma(s)-\sigma(t)||\leq\int_{s}^{t}\parallel\dot{\sigma}(\tau)\parallel_{L^{2}(I)}d\tau,$$
we conclude $$\sigma\in AC^{2}(0,T;L^{2}(I)).$$
\begin{flushright}
  $\Box$
\end{flushright}

Next, we will prove the Tonelli Theorem by the method of direct variation. As we know that the method needs that the action functional is lower semicontinuious and that the variation space is compact with respect to some topology. Firstly, we define the action functional
\begin{equation}\label{e31}
  \mathbb{L}(\sigma) :=\int_{0}^{T}L(\sigma(t),\dot{\sigma}(t))dt=\int_{0}^{T}\frac{1}{2}\parallel \dot{\sigma}(t)\parallel_{L^{2}(I)}^{2}-\mathcal{W}(\sigma(t))dt.
\end{equation}
\begin{proposition}   \label{p31}
 Let $\sigma^{k},\sigma\in AC^{2}(0,T;L^{2}(I)).$ If $\sigma^{k}\stackrel{\tau}{\rightarrow}\sigma,$ then
\begin{equation*}
  \liminf_{k} \mathbb{L}(\sigma^{k}) \geq  \mathbb{L}(\sigma).
\end{equation*}
\end{proposition}
\noindent\textbf{Proof.} The first step is to prove that
\begin{equation}\label{e32}
   \liminf_{k}\int_{0}^{T}\parallel \dot{\sigma}^{k}(t)\parallel_{L^{2}(I)}^{2}dt\geq\int_{0}^{T}\parallel \dot{\sigma}(t)\parallel_{L^{2}(I)}^{2}dt.
 \end{equation}
Let $l:=\liminf_{k}\int_{0}^{T}\parallel \dot{\sigma}^{k}(t)\parallel_{L^{2}(I)}^{2}dt\geq0$. If $l=+\infty$, (\ref{e32}) is obvious, otherwise we can choose a subsequence still denoted by $\sigma^{k}$ such that
$$\lim_{k}\int_{0}^{T}\parallel \dot{\sigma}^{k}(t)\parallel_{L^{2}(I)}^{2}dt=l\quad\textrm{and}\quad\int_{0}^{T}\parallel \dot{\sigma}^{k}(t)\parallel_{L^{2}(I)}^{2}dt\leq(\sqrt{l}+\frac{1}{k})^{2}.$$
\\For any $\varphi\in C_{c}^{\infty}([0,T]\times I)$, we have
\begin{eqnarray*}
 <\dot{\sigma},\varphi>_{L^{2}([0,T]\times I)} &=& -<\sigma,\dot{\varphi}>_{L^{2}([0,T]\times I)}\\
 &=& -\int_{0}^{T}<\sigma(t),\dot{\varphi}(t)>_{L^{2}(I)}dt \\
   &=&\lim_{k}-\int_{0}^{T}<\sigma^{k}(t),\dot{\varphi}(t)>_{L^{2}(I)}dt \\
  &=&\lim_{k}\int_{0}^{T}<\dot{\sigma}^{k}(t),\varphi(t)>_{L^{2}(I)}dt \\
   &\leq& \lim_{k}(\int_{0}^{T}\parallel \dot{\sigma}^{k}(t)\parallel_{L^{2}(I)}^{2}dt)^{\frac{1}{2}}\parallel\varphi\parallel_{L^{2}([0,T]\times I)} \\
 &\leq& \lim_{k}(\sqrt{l}+\frac{1}{k})\parallel\varphi\parallel_{L^{2}([0,T]\times I)},
\end{eqnarray*}
where the third line has used the Lebesgue Dominated Convergence Theorem. Hence
$$\int_{0}^{T}\parallel \dot{\sigma}(t)\parallel_{L^{2}(I)}^{2}dt=\parallel\dot{\sigma}\parallel_{L^{2}([0,T]\times I)}^{2}\leq l.$$
\begin{flushright}
  $\Box$
\end{flushright}
Secondly, we need to prove
\begin{equation}\label{e33}
  \limsup_{k}\int_{0}^{T}\mathcal{W}(\sigma^{k}(t)) dt\leq\int_{0}^{T}\mathcal{W}(\sigma(t))dt.
\end{equation}
Since $\mathcal{W}$\;is bounded above on $L^{2}(I)$, it is easy to prove (\ref{e33}) true by Fatou Lemma, refer to \cite{gan3}.
\begin{flushright}
  $\Box$
\end{flushright}

Up to now, Proposition \ref{p31} and Lemma \ref{l32} have paved the way for the proof of Tonelli Theorem.
\begin{theorem}\textbf{(Tonelli Theorem)}   \label{t31}
 For any $M,N \in L^{2}(I)$, let $AC^{2}_{M,N}(0,T;L^{2}(I)):=\{\sigma\in AC^{2}(0,T;L^{2}(I))\;\mid\sigma_{0}=M,\sigma_{T}=N \}$. Then $\mathbb{L}(\sigma)$ has a minimizer in $ AC^{2}_{M,N}(0,T;L^{2}(I))$ .
\end{theorem}
\noindent\textbf{Proof.} Set
$$C_{inf}=\inf_{\sigma\in AC^{2}_{M,N}(0,T;L^{2}(I))}\mathbb{L}(\sigma).$$
We know that $C_{inf}<+\infty$, indeed, choose the line $\sigma^{*}$ between $M$ to $N$, then
\begin{equation*}
   C_{inf} \leq \mathbb{L}(\sigma^{*}) \leq \frac{||M-N||^{2}_{L^{2}(I)}}{2T}+K_{0}T.
\end{equation*}
Consider the subset $\Sigma_{1}$ of $AC^{2}_{M,N}(0,T;L^{2}(I))$ formed by the curves $\sigma$ such that $\mathbb{L}(\sigma)\leq C_{inf}+1.$ This set is nonempty by definition. It is also closed according to the Proposition \ref{p31}. Due to the boundedness of the potential $\mathcal{W}$, there exists $C$ such that
 $$\int_{0}^{T}\parallel \dot{\sigma}(t)\parallel_{L^{2}(I)}^{2}dt\leq C,\; \forall \sigma\in \Sigma_{1}.$$
Therefore, $\Sigma_{1}$ is nonempty compact with respect to topology $\tau$ by Lemma \ref{l32}. And together with the lower semicontinuity of $\mathbb{L}(\sigma)$, we conclude the theorem by the method of direct variation.
\begin{flushright}
  $\Box$
\end{flushright}

\section{The Lax-Oleinik semigroup}
We introduce a semigroup of nonlinear operators $\{T_{t}^{-}\}_{t\geq0}$. Given a function $u:L^{2}(I)\rightarrow[-\infty,+\infty]$ and $t>0$, we define a function
$$T_{t}^{-}u:L^{2}(I)\rightarrow[-\infty,+\infty]$$
by
\begin{equation}\label{e41}
 T_{t}^{-}u(M)=\inf_{\sigma}\{u(\sigma(0))+\int_{0}^{t}L(\sigma(s),\dot{\sigma}(s))ds\},
\end{equation}
where the infimum is taken on all 2-absolutely continuous curves $\sigma:[0,t]\rightarrow L^{2}(I)$ with $\sigma(t)=M.$ In the following lemma, we show some properties of the nonlinear operators.
\begin{lemma}\label{l42}
\begin{description}
  \item[(1)]  We have $T_{t+s}^{-}=T_{t}^{-}\circ T_{s}^{-},$ \;for each $t,s\geq0$.
  \item[(2)]  For every $u,v:L^{2}(I)\rightarrow[-\infty,+\infty]$ and all $t\geq 0,$ one gets that $T_{t}^{-}u\leq T_{t}^{-}v$ if $u\leq v$.
  \item[(3)]  If $k\in\mathbb{R}$ and $u:L^{2}(I)\rightarrow[-\infty,+\infty]$, we have $T_{t}^{-}(u+k)=T_{t}^{-}(u)+k$.
\end{description}
\end{lemma}
The assertions (1),(2) and (3) are not difficult to prove by the definition of $T_{t}^{-}$. From this lemma, we know that the operators is a semigroup which is well known in PDE and in Calculus of Variations. People call it Lax-Oleinik semigroup. In the following part, we further discuss the semigroup.

\begin{proposition}\label{p41}
If $t\geq0$ and $u:L^{2}(I)\rightarrow[-\infty,+\infty]$ is periodic and rearrangement invariant, then so is $T^{-}_{t}u.$
\end{proposition}
The proposition is easy to check by definition.
\begin{proposition}\label{p42}
The maps $T^{-}_{t}$ are non-expansive. Namely, if $u,v$ are bounded functions defined on $L^{2}(I)$, then $||T^{-}_{t}u-T^{-}_{t}v||_{\infty}\leq ||u-v||_{\infty}$, where $||\cdot||_{\infty}$ is the $L^{\infty}$ norm.
\end{proposition}
\noindent\textbf{Proof.}
Notice that $-||u-v||_{\infty}+v\leq u\leq ||u-v||_{\infty}+v$, and we use the assertions (2) and (3) in Lemma \ref{l42} to find
\begin{equation*}
  -||u-v||_{\infty}+T^{-}_{t}v\leq T^{-}_{t}u\leq ||u-v||_{\infty}+T^{-}_{t}v.
\end{equation*}
Thus, we conclude the proposition.
\begin{flushright}
  $\Box$
\end{flushright}
\begin{proposition}\label{p43}
Assume that a bounded function $u:L^{2}(I)\rightarrow[-\infty,+\infty]$ is weakly lower semicontinuous, $T^{-}_{t}u$ is also weakly lower semicontinuous.
\end{proposition}
\noindent\textbf{Proof.}
Suppose that a sequence $M_{n}$ in $L^{2}(I)$ weakly converge to $M$, we need to prove that
\begin{equation*}
  \liminf_{n} T_{t}^{-}u(M_{n})\geq T_{t}^{-}u(M).
\end{equation*}
Set $l:=\liminf_{n} T_{t}^{-}u(M_{n})$, obviously, $l<+\infty.$ One can choose a subsequence still denoted by $M_{n}$ such that $\lim_{n} T_{t}^{-}u(M_{n})=l$ and $T_{t}^{-}u(M_{n})\leq l+1,\;\forall n\geq 1$. For each $ T_{t}^{-}u(M_{n})$, by definition, there exists a curve $\sigma^{n}\in AC^{2}(0,t;L^{2}(I))$ with $\sigma^{n}(t)=M_{n}$ in such way that
\begin{equation}\label{e42}
 T_{t}^{-}u(M_{n})>u(\sigma^{n}(0))+\int_{0}^{t}L(\sigma^{n}(s),\dot{\sigma}^{n}(s))ds-\frac{1}{n}.
\end{equation}
Recall the boundedness of $u$ and the potential $\mathcal{W}$, we find a constant $C$ such that
 $$\int_{0}^{t}\parallel \dot{\sigma}^{n}(s)\parallel_{L^{2}(I)}^{2}ds\leq C,\; \forall n\geq 1.$$
Since $\sigma^{n}(t)=M_{n}\rightharpoonup M,$ by Lemma \ref{l32}, there exist $\sigma\in AC^{2}(0,t;L^{2}(I))$ and a subsequence still denoted by $\sigma^{n}$ such that $\sigma^{n}\stackrel{\tau}{\rightarrow}\sigma.$ It is clearly know $\sigma(t)=M.$ Letting $n\rightarrow+\infty$ in (\ref{e42}), and combining with Proposition \ref{p31}, we obtain
\begin{equation*}
  l\geq u(\sigma(0))+\int_{0}^{t}L(\sigma(s),\dot{\sigma}(s))ds\geq T_{t}^{-}u(M).
\end{equation*}
\begin{flushright}
  $\Box$
\end{flushright}
\begin{proposition}\label{p44}
$u$ is under the same assumption as that in above proposition, if $t>0$ and $M\in L^{2}(I)$ are given, there exists $\sigma\in AC^{2}(0,t;L^{2}(I))$ such that $\sigma(t)=M$ and $$ T_{t}^{-}u(M)=u(\sigma(0))+\int_{0}^{t}L(\sigma(s),\dot{\sigma}(s))ds$$
\end{proposition}

The proof is similar to that of Proposition \ref{p43}. We only replace the left part of (\ref{e42}) by $ T_{t}^{-}u(M)$. The rest part is almost the same.

\section{The dominated functions}
If $u$ is a weak KAM solution for Lagrangian function $L$ defined on $L^{2}(I)$, then $u$ is dominated by $L+c$, denoted by $u\prec L+c$, that is, for each $\sigma \in AC^{2}(a,b;L^{2}(I))$, we have
\begin{equation}\label{e51}
  u(\sigma(b))-u(\sigma(a))\leq \int_{a}^{b}L(\sigma(s),\dot{\sigma}(s))ds+c(b-a).
\end{equation}
Remark that $a,b$ in above inequality is also arbitrary. Given $c\in\mathbb{R}^{1}$, set $\mathcal{H}(c):=\{u: L^{2}(I)\rightarrow\mathbb{R}^{1}\;|\;u\prec L+c\}.$ We will search a weak KAM solution in $\mathcal{H}(c)$. Firstly, we display some properties of $\mathcal{H}(c)$.
\begin{lemma}
\begin{enumerate}
  \item Fix any $k\in\mathbb{R}$, then $u\in\mathcal{H}(c)$ if and only if $u+k\in\mathcal{H}(c)$;
  \item  Every $u$ in $\mathcal{H}(c)$ is $(\frac{1}{2}+K_{0}+c)-$Lipschitz continuous, namely,
             $$|u(\bar{M})-u(M)|\leq (\frac{1}{2}+K_{0}+c)||\bar{M}-M||_{L^{2}(I)};$$
  \item If $u$ is K-Lipschitz continuous, then $u\in\mathcal{H}(\frac{K^{2}}{2}+K_{0})$;
  \item The set $\mathcal{H}(c)$ is closed convex for the topology of pointwise convergence.
\end{enumerate}
\end{lemma}
\noindent\textbf{Proof.}
Statement 1 and 4 are easy to check by definition. \\
To show statement 2, we choose the line $\tilde{\sigma}:[0,||\bar{M}-M||_{L^{2}(I)}]\rightarrow L^{2}(I) $ parametrized by unit length such that $\tilde{\sigma}(0)=M,$ $\tilde{\sigma}(||\bar{M}-M||_{L^{2}(I)})=\bar{M}.$ Due to (\ref{e51}), we obtain
\begin{eqnarray*}
    u(\tilde{\sigma}(||\bar{M}-M||_{L^{2}(I)}))-u(\tilde{\sigma}(0))&\leq& \int_{0}^{||\bar{M}-M||_{L^{2}(I)}}L(\tilde{\sigma}(s),\dot{\tilde{\sigma}}(s))ds+c||\bar{M}-M||_{L^{2}(I)} \\
  &\leq &   (\frac{1}{2}+K_{0}+c)||\bar{M}-M||_{L^{2}(I)}.
\end{eqnarray*}
Exchanging the role of $M$ and $\bar{M}$, we conclude the statement 2.\\
As to statement 3, for any $\sigma \in AC^{2}(a,b;L^{2}(I))$, we know
\begin{eqnarray*}
  u(\sigma(b))-u(\sigma(a)) &\leq& K||\sigma(b)-\sigma(a)||_{L^{2}(I)} \\
   &\leq & \int_{a}^{b}K\parallel \dot{\sigma}(s)\parallel_{L^{2}(I)}ds  \\
  &\leq &  \frac{1}{2}\int_{a}^{b}K^{2}+\parallel \dot{\sigma}(s)\parallel_{L^{2}(I)}^{2}ds\\
   &\leq &  \int_{a}^{b}L(\sigma(s),\dot{\sigma}(s))ds+(\frac{K^{2}}{2}+K_{0})(b-a).
\end{eqnarray*}
Therefore, $u\in \mathcal{H}(\frac{K^{2}}{2}+K_{0}).$
\begin{flushright}
  $\Box$
\end{flushright}
\begin{lemma}\label{l52}
\begin{enumerate}
  \item For all $t\geq 0,$ $u\prec L+c$ if and only if $u\leq T^{-}_{t}u+ct$.
  \item The maps $T^{-}_{t}$ send $\mathcal{H}(c)$ into itself.
  \item Assume that $u\in\mathcal{H}(c)$ is bonuded, then the map $t\rightarrow T^{-}_{t}u$ is continuous on $[0,+\infty)$, where we choose $C^{0}$ topology in $\mathcal{H}(c)$.
 \end{enumerate}
\end{lemma}
\noindent\textbf{Proof.}
The first statement is immediate from the definitions. Using Lemma \ref{l42} and former statements it is not difficult to verify the second one. It remains to prove the last one. Since $T^{-}_{t}$ is a semigroup, it is enough to show that the map is continuous at $t=0$. For all $t> 0$ and $M\in L^{2}(I),$ by the first statement, it is clear that
\begin{equation*}
   T^{-}_{t}u(M)\geq u(M)-ct.
\end{equation*}
Conversely, using the constant curve $\sigma:[0,t]\rightarrow L^{2}(I), s\mapsto M$, we obtain
\begin{equation*}
T_{t}^{-}u(M)\leq u(M) +K_{0}t.
\end{equation*}
So we have
\begin{equation*}
  ||T^{-}_{t}u-u||_{\infty}\leq t\max\{|c|,K_{0}\}.
\end{equation*}
\begin{flushright}
  $\Box$
\end{flushright}
\section{The proof of weak KAM Theorem}
This section is mainly divided into two parts. One part is to prove that the operators $T^{-}_{t}$ have a common fixed point. The other part is devoted to checking that the common fixed point is a weak KAM solution. Firstly, we recall a fixed point theorem.
\begin{lemma}\label{l61}
Let $E$ be a closed convex space and $\varphi_{t}:E\rightarrow E$ be a family of maps defined for $t\geq 0.$ We assume that the following conditions are satisfied
\begin{enumerate}
  \item For each $t,s\geq 0,$ we have $\varphi_{t+s}=\varphi_{t}\circ\varphi_{s}.$
  \item For each $t\geq 0,$ the map $\varphi_{t}$ is non-expansive.
  \item Fixed $x\in E$, the map $t\rightarrow \varphi_{t}(x)$ is continuous on $[0,+\infty).$
  \item For each $t> 0,$ the image $\varphi_{t}(E)$ is relatively compact.
\end{enumerate}
Then the maps $\varphi_{t}$ have a common fixed point.
\end{lemma}
To apply this lemma, we construct a suitable space
\begin{eqnarray*}
  E: &=& \mathcal{H}(c)\bigcap\{u: L^{2}(I)\rightarrow\mathbb{R}\;|\; u \text{ is periodic, rearrangement invariant, }\\
   & &\quad\quad\quad\quad\quad\quad\text{ and weakly lower semicontinuous.} \}.
\end{eqnarray*}
The norm on $E$ is $C^{0}$ norm. For $c$ large enough, $E$ is nonempty. By Remarks \ref{r21} and \ref{r22}, we know that $\forall u\in E$, $||u||_{\infty}<\infty$.
\begin{lemma}
$E$ constructed above is closed convex.
\end{lemma}
\noindent\textbf{Proof.}
It is sufficient to check that $E$ is closed. We suppose that a sequence $u_{n}$ belongs to $E$ and $u_{n}\stackrel{C^{0}}{\rightarrow} u$. It is obvious that $u$ is periodic, rearrangement invariant. So we just check that $u$ is weakly lower semicontinuous. That is, if $M_{m}\rightharpoonup M_{0}$, we prove
\begin{equation}\label{e61}
 \liminf_{m}u(M_{m})\geq u(M_{0}).
\end{equation}
$\forall \epsilon>0,$ since $u_{n}\stackrel{C^{0}}{\rightarrow} u$, there exists $n_{0}$ such that if $n>n_{0},$
\begin{equation*}
  |u_{n}(M)-u(M)|<\epsilon,\quad\forall M\in L^{2}(I).
\end{equation*}
Fix some $n>n_{0}$, according to the assumption that $u_{n}$ is weakly lower semicontinuous, there exists $m_{0}$ such that if $m>m_{0},$
\begin{equation*}
  u_{n}(M_{m})>u_{n}(M_{0})-\epsilon.
\end{equation*}
Therefore, if $m>m_{0}$ and $n>n_{0},$ one gets
\begin{eqnarray*}
  u(M_{m})-u(M_{0}) &=& u(M_{m})-u_{n}(M_{m})+ u_{n}(M_{m})-u_{n}(M_{0})+ u_{n}(M_{0})-u(M_{0}) \\
   &>& -3\epsilon.
\end{eqnarray*}
Then letting $m\rightarrow \infty,$ to find
\begin{equation*}
  \liminf_{m}u(M_{m})\geq u(M_{0})-3\epsilon
\end{equation*}
which induces (\ref{e61}).
\begin{flushright}
  $\Box$
\end{flushright}
\begin{theorem}
If the potential $\mathcal{W}$ is continuous, periodic, rearrangement invariant and weakly lower semicontinuous, then there exist $u\in E$ and a constant $\lambda$ such that $u=T^{-}_{t}u+\lambda t$, $\forall t\geq 0.$
\end{theorem}
\noindent\textbf{Proof.}
By Proposition \ref{p41}, \ref{p42} and Lemma \ref{l52}, the operators $T^{-}_{t}$ map $E$ into itself. Using Lemma \ref{l42}, Proposition \ref{p41} and Lemma \ref{l52}, it is obvious that $T^{-}_{t}$ satisfy the first three conditions in Lemma \ref{l61}, however, the last one is not satisfied. To overcome this difficulty, we consider the quotient $\bar{E}:=E/\mathbb{R}.$ This quotient space $\bar{E}$ is also convex closed for the quotient norm
\begin{equation*}
  ||[u]||=\inf_{a\in\mathbb{R}}||u+a||_{\infty},
\end{equation*}
where $[u]$ is the class in $\bar{E}$ of $u\in E.$ We introduce $E_{0}$ the subset of $E$ whose elements vanish at $0\in L^{2}(I).$ By Remarks \ref{r21} and \ref{r22},  $\forall u\in E_{0}$,
\begin{equation*}
 ||u||_{\infty}=\max_{M\in L^{2}(I)}|u(M)|=\max_{M\in L^{2}(I)}|u(M)-u(0)|\leq (\frac{1}{2}+K_{0}+c)Diam (\mathcal{P}(\mathbb{T}^{1})).
\end{equation*}
Since the functions in $E_{0}$ are equi-Lipschitz, it follows from Ascoli-Arzal\`{a} Theorem that $E_{0}$ is compact. Besides, it is easy to check $E_{0}/\mathbb{R}=\bar{E}$. Thus, $\bar{E}$ is compact for the quotient topology. Recall that $T_{t}^{-}(u+a)=T_{t}^{-}(u)+a$, we know that the maps $T_{t}^{-}$ pass the quotient to a semigroup $\bar{T}_{t}^{-}:\bar{E}\rightarrow\bar{E}$. Now the image of $\bar{T}_{t}^{-}$ is compact. Using Lemma \ref{l61}, we find a common fixed point for all $\bar{T}_{t}^{-}$. We then deduce that there exists $u\in E$ such that $u=T^{-}_{t}u+c_{t},$ where $c_{t}$ is a constant. Due to the semigroup property, we get $c_{t+s}=c_{t}+c_{s}$. Moreover, $c_{t}=\lambda t$ with $\lambda=c_{1}$ because of the continuity of the map $t\rightarrow T_{t}^{-}u$.
\begin{flushright}
  $\Box$
\end{flushright}
\begin{theorem}\label{t64}
If $u=T^{-}_{t}u+\lambda t,\;\forall t\geq 0$, then $u\prec L+\lambda$ and for each $M\in L^{2}(I),$ there exists a curve $\sigma^{M}:(-\infty,0] \rightarrow L^{2}(I)$ with $\sigma^{M}(0)=M$ such that
\begin{equation*}
  u(M)-u(\sigma^{M}(-t))=\int_{-t}^{0}L(\sigma^{M}(s),\dot{\sigma}^{M}(s))ds+\lambda t, \;t\geq 0.
\end{equation*}
\end{theorem}
\noindent\textbf{Proof.}
By Lemma \ref{l52}, we easily get $u\prec L+\lambda$. It remains to show the existence of calibrated curve for any given $M\in L^{2}(I).$ We already know that, for each $t>0$, there exists a curve $\bar{\sigma}_{t}:[0,t]\rightarrow L^{2}(I),$ with $\bar{\sigma}_{t}(t)=M$ and
\begin{equation}\label{e62}
 u(M)-\lambda t= T_{t}^{-}u(M)=u(\bar{\sigma}_{t}(0))+\int_{0}^{t}L(\bar{\sigma}_{t}(s),\dot{\bar{\sigma}}_{t}(s))ds.
\end{equation}
Set $\sigma_{t}(s)=\bar{\sigma}_{t}(s+t),$ then (\ref{e62}) turns into
\begin{equation*}
 u(M)-u(\sigma_{t}(-t))=\int_{-t}^{0}L(\sigma_{t}(s),\dot{\sigma}_{t}(s))ds+\lambda t, \;\forall t\geq 0.
\end{equation*}
Moreover, we know that
\begin{equation}\label{e63}
 u(M)-u(\sigma_{t}(-t'))=\int_{-t'}^{0}L(\sigma_{t}(s),\dot{\sigma}_{t}(s))ds+\lambda t', \;\forall t'\in [0,t].
\end{equation}
Especially, for each positive integer $n,$ we obtain a curve $\sigma_{n}:[-n,0]\rightarrow L^{2}(I),$ with $\sigma_{n}(0)=M$ and satisfies (\ref{e63}). Fix $n>0$, when $m>n$, replacing $t,t'$ by $m,n$ in (\ref{e63}), to find
\begin{equation}\label{e64}
 u(M)-u(\sigma_{m}(-n))=\int_{-n}^{0}L(\sigma_{m}(s),\dot{\sigma}_{m}(s))ds+\lambda n.
\end{equation}
Since $u$ is bounded, there exists a constant $C$ (depends on $n$ and $||u||_{\infty}$) such that
 \begin{equation*}
  \int_{-n}^{0}\parallel \dot{\sigma}_{m}(s)\parallel_{L^{2}(I)}^{2}ds\leq C.
 \end{equation*}
Then, following Lemma \ref{l32}, the curves $\sigma_{m}|_{[-n,0]}$ are compact with respect to the topology $\tau$ in Definition \ref{d31}. By a diagonal process, we can extract an increasing sequence of indexes $m_{k}\in\mathbb{N}$ such that, for each $n>0$, the sequence $(\sigma_{m_{k}}|_{[-n,0]})_{m_{k}>n}$ converges when $k\rightarrow\infty$. Define a curve $\sigma^{M}:(-\infty,0]\rightarrow L^{2}(I)$ by $\sigma(t)=\lim_{k}\sigma_{m_{k}}(t)$. Observe now that, each curve $(\sigma_{m_{k}}|_{[-n,0]})_{m_{k}>n}$ satisfies (\ref{e64}). The lower semicontinuity of the action functional and function $u$ implies
\begin{equation*}
 u(M)-u(\sigma^{M}(-t))\geq\int_{-t}^{0}L(\sigma^{M}(s),\dot{\sigma}^{M}(s))ds+\lambda t, \;t\geq 0.
\end{equation*}
Therefore, we obtain the conclusion by the fact that $u\prec L+\lambda$.

\section*{Acknowledgements}

The authors are very grateful to Professor Xiaoping Yuan and Xuemei Li for their invaluable discussions
and encouragements.


\begin{thebibliography}{aa}



\bibitem{amb}L. Ambrosio, N. Gigli and G. Savar\'{e}. Gradient flows in metric spaces and the Wasserstein spaces
                     of probability measures. Lectures in Mathematics, ETH Zurich, Birkh\"{a}user, (2005).
          \bibitem{con} G. Contreras, R. Iturriaga, G. Paternain, and M. Paternain. Lagrangian graphs, minimizing measures and Ma$\tilde{\text{n}}$\'e critical values. Geom. Funct. Anal. 8 (1998), 788-809.
          \bibitem{cra} M. G. Crandall and P. L. Lions. Viscosity solutions of Hamilton-Jacobi equations. Trans. Amer. Math. Soc. 277 (1983), l-42.
          \bibitem{fat1} A. Fathi. Solutions KAM faibles conjuguees et barrieres de Peierls. C. R. Acad. Sci. Paris Ser.I Math. 325 (1997), 649-652.
          \bibitem{fat2} A. Fathi. Theoreme KAM faible et theorie de Mather sur les systemes lagrangiens. C. R. Acad.Sci. Paris Ser. I Math. 324 (1997), 1043-1046.
          \bibitem{fat3} A. Fathi. Sur la convergence du semi-groupe de Lax-Oleinik. C. R. Acad. Sci. Paris Ser. I Math.327 (1998), 267-270.
          \bibitem{fat4} A. Fathi and E. Maderna.Weak KAM theorem on non compact manifolds. NoDEA Nonlinear Differential Equations Appl. 14 (2007), 1-27.
          \bibitem{fat5} A. Fathi. Weak KAM theorem and Lagrangian dynamics. 10th Preliminary version, (2009).
          \bibitem{fat6} A. Fathi and John N. Mather. Failure of the convergence of the Lax-Oleinik semigroup in the time periodic case. Bull. Soc. Math. France 128 (2000), 473-483.
          \bibitem{gan1} W. Gangbo, A. Tudorascu. Lagrangian Dynamics on an infinite-dimensional torus: a Weak KAM theorem. Adv. Math. 224 (2010), 260-292.
          \bibitem{gan2} W. Gangbo, A. Tudorascu. Weak KAM on the wasserstein torus with multi-dimensional underlying space. Comm. Pure Appl. Math. 67 (2014), 408-463.
          \bibitem{gan3} W. Gangbo, N. Tguyen  and  A. Tudorascu. Hamilton-Jacobi equations in the Wasserstein space. Methods and Applications Analysis 15 (2008), 155-184.
          \bibitem{mad} E. Maderna. On weak KAM theory for N-body problems. Ergodic Theory Dynam. Systems 32 (2012), 1019¨C1041.
          \bibitem{wan} K.Wang and J. Yan. A new kind of Lax-Oleinik type operator with parameters for time-periodic positive definite Lagrangian systems. Commun. Math. Phys. 309 (2012), 663-691.





\end{thebibliography}
  \end{document}